\documentclass[12pt]{article}
\usepackage{graphicx,psfrag,epsfig}
\usepackage{amssymb,amsmath,amscd,amsthm}
\usepackage{graphicx,psfrag,epsfig}

\usepackage{graphicx}
\usepackage[active]{srcltx}

\newtheorem{theorem}{Theorem}[section]

\newtheorem{lemma}[theorem]{Lemma}

\date{}

\begin{document}

\date{}
\title{Branching
diffusion in the super-critical regime}
\author{L.
Koralov\footnote{Dept of Mathematics, University of Maryland,
College Park, MD 20742, koralov@math.umd.edu} , S.
Molchanov\footnote{Dept of Mathematics, University of North
Carolina, Charlotte, NC 28223, smolchan@uncc.edu} } \maketitle
\begin{abstract}
We investigate the long-time evolution of branching diffusion
processes (starting with a single particle) in inhomogeneous
media. The qualitative behavior of the processes depends on the
intensity of the branching. We analyze the super-critical case,
when the total number of particles growing exponentially with
positive probability. We study the asymptotics of the number of
particles in different regions of space and describe the growth of
the region occupied by the particles.
\end{abstract}

%82B26,
%82B27, 82D60, 35K10.
\section{Introduction}

The mathematical study of branching processes goes back to the
work of Galton and Watson \cite{WG} who were interested in the
probabilities of long-term survival of family names. Later it was
realized that similar mathematical models could be used to
describe the evolution of a variety of biological populations, in
genetics \cite{Fi1, Fi2, Fi3, Hal}, and in the study of certain
chemical and nuclear reactions \cite{Se, HaU}. The branching
processes (in particular, branching diffusions) are central in the
study of the evolution of various populations such as bacteria,
cancer cells, carriers of a particular gene, etc., where each
member of the population may die or produce offspring
independently of the rest.

In this paper we describe the long-time behavior of the population
in different regions of space when, in addition to branching, the
members of the population move diffusively in space and the
branching mechanism depends on the location. In particular, we'll
consider regions $U_{x(t)}$ centered at a point $x(t) = v t$ which
is at a linear (in $t$) distance from the origin at time $t$. We
will be interested in the super-critical case (when the total
population grows exponentially with positive probability).

 Consider a collection
of particles in $ \mathbb{R}^d$ that move diffusively and
independently starting with one particle. Besides the diffusive
motion, the particles can duplicate with the rate of duplication
$v(x)$, $x \in \mathbb{R}^d$, where $x$ is the position of a given
particle and $v$ is a continuous non-negative compactly supported
function. Both copies start moving independently immediately after
the duplication (annihilation of particles and creation of more
than two particles from one could also be considered, but here we
discuss only duplication for clarity of exposition).

Let $n^x_t(U)$ be the number of particles in a domain $U \subseteq
\mathbb{R}^d$ at time $t$, assuming that at time zero there was a
single particle located at $x$. The large time behavior of
$n^x_t(U)$ depends crucially on the magnitude of $v$, that is on
whether the operator
\begin{equation} \label{opera}
\mathcal{ L} u(x) = \frac{1}{2} \Delta u(x) + v(x) u(x)
\end{equation}
has a positive eigenvalue. This is the operator in the right hand
side of the equations on the particle density and higher order
correlation functions, given below. In the super-critical case
(i.e., if there exists a positive eigenvalue),  if $U$ is fixed,
the results of \cite{EHK} on the asymptotics of $n^x_t(U)$ cover,
in particular, the case of compactly supported $v$. Namely,
$n^x_t(U)$ grows exponentially with a random coefficient in front
of the exponent. The random coefficients corresponding to
different domains differ only by a multiplicative constant. (See
also \cite{W}, \cite{AH}). The distribution of the random
coefficient can be described in terms of its moments (see
\cite{K}). (For the asymptotic properties of branching random
walks see also \cite{ABY}, \cite{BY1}, \cite{ABY2}, \cite{Ya}.) If
$\mathcal{L}$ has no positive eigenvalues, the total number of
particles tends to a finite random limit whose distribution can be
described in terms of its moments.

In the current paper we consider the super-critical case. We use
the spectral techniques developed in \cite{CKMV1}, \cite{K} to get
the asymptotic formulas for the density and higher order
correlation functions of the branching process. These allow us to
get the asymptotics of $n^x_t(U_{tv})$, where $U_{tv}$ is a region
of fixed size centered at a point whose distance from the origin
grows linearly with $t$. This asymptotics is the main result of
this paper.

After recalling the equations on the correlation functions
(Section \ref{momeq}) and the asymptotic behavior of the
correlation functions (Section \ref{sectt}), we show in Section
\ref{growth} that the total number of particles, after division by
an exponential factor, tends to a random limit in $L^2$. In
Sections~\ref{growth22} and \ref{asco} we prove a similar result
for fixed domains and for domains located at a linear (in $t$)
distance from the origin and show that the convergence takes place
not only in $L^2$ but also almost surely. In Section \ref{shape}
we show that on the event that the number of particles grows
exponentially, the region occupied by the particles grows linearly
in $t$. In Section \ref{finitely} we give the distribution of the
limiting number of particles in the event that the limiting number
of particles is finite.

\section{Equations on correlation functions}
\label{momeq}

%Since we deal in this paper with the case when $\beta > \beta_{\rm
%cr}$, where $\beta$ is kept fixed, we'll incorporate the constant
%$\beta$ into the function $v$. In particular, instead of
%(\ref{opera}) we'll write
%\begin{equation} \label{opera2}
% L u(x) = \frac{1}{2} \Delta u(x) + v(x) u(x)
%\end{equation}
%and assume that $L$ has a positive eigenvalue.
%
 Let
$B_\delta$ be a ball of radius $\delta$ in $ \mathbb{R}^d$. For $t
> 0$ and $x,y_1,y_2,... \in \mathbb{R}^d$ with all $y_i$  distinct,
define the particle density $\rho_1(t,x,y_1)$ and the higher order
correlation functions $\rho_n(t,x,y_1,...,y_n)$ as the limits of
probabilities of finding $n$ distinct particles in
$B_\delta(y_1)$,...,$B_\delta(y_n)$, respectively, divided by
${\rm Vol}^n (B_\delta)$, under the condition that there is a
unique particle at $t = 0$ located at $x$. We extend
$\rho_n(t,x,y_1,...,y_n)$ by continuity to allow for $y_i$ which
are not necessarily distinct. For fixed $y_1$, the density
satisfies the equation
\begin{equation} \label{fmo}
\partial_t \rho_1(t,x,y_1) = \frac{1}{2} \Delta \rho_1(t,x,y_1) + v(x) \rho_1(t,x,y_1),
\end{equation}
\[
\rho_1(0,x,y_1) = \delta_{y_1}(x).
\]
Indeed, let $s,t > 0$. Then we can write
\begin{equation} \label{appro}
\rho_1(s+t,x,y_1) = (2 \pi s)^{-\frac{d}{2}} \int_{ \mathbb{R}^d}
e^{\frac{|x - z|^2}{2s}} \rho_1(t,z,y_1)d z +  v(x) s
\rho_1(t,x,y_1) + \alpha(s,t,x,y_1),
\end{equation}
 where the term with the
integral on the right hand side is due to the effect of the
diffusion on the interval $[0,s]$, the second term is due to the
probability of branching on $[0,s]$, and $\alpha$ is the
correction term. The correction term is present since (a) more
than one instance of branching may occur before time $s$, and (b)
even if a single branching occurs between the times $0$ and $s$,
then the original particle will be located not at $x$ but at a
nearby point and the intensity of branching there is slightly
different from $v(x)$. It is clear that $\lim_{s \downarrow 0}
\sup_{x,y \in \mathbb{R}^d} \alpha(s,t,x,y)/s = 0$. After
subtracting $\rho_1(t,x,y_1)$ from both sides of~(\ref{appro}),
dividing by $s$ and taking the limit as $s \downarrow 0$, we
obtain (\ref{fmo}).

%This can be derived by examining the behavior of the original
%particle on the time interval $[0,\Delta t]$, as $\Delta t
%\downarrow 0$, with the first term on the right hand side coming
%from the effect of the diffusion and the second term from the
%branching.
The equations on $\rho_n$, $n > 1$, are somewhat more complicated:
%\begin{equation} \label{fmo2}
\begin{equation} \label{manyp}
\partial_t \rho_n(t,x,y_1,...,y_n) = \frac{1}{2} \Delta \rho_n(t,x,y_1,...,y_n) +
 v(x) \left( \rho_n(t,x,y_1,...,y_n) +  H_n (t,x,y_1,...,y_n)
\right) ,
\end{equation}
%\end{equation}
%\begin{equation} \label{ico2}
\[
\rho_n(0,x,y_1,...,y_n) \equiv 0.
\]
%\end{equation}
Here
\[
H_n (t,x,y_1,...,y_n) = \sum_{U \subset Y, U \neq \emptyset}
\rho_{|U|} (t,x,U) \rho_{n-|U|}(t,x, Y\setminus U),
\]
where $Y = (y_1,...,y_n)$, $U$ is a proper non-empty subsequence
of $Y$, and $|U|$ is the number of elements in this subsequence.
Equation (\ref{manyp}) is derived similarly to (\ref{fmo}). The
combinatorial term $H_n$ appears after taking into account the
event that there is a single branching on the time interval
$[0,s]$, the descendants of the first particle are found at the
points in $U$ at time $s+t$, while the descendants of the second
particle are found at the points of $ Y \setminus U$, with the
summation over all possible choices of $U$.

\section{Asymptotics of the correlation functions} \label{sectt}

First we recall some basic facts about the operator
$\mathcal{L}:L^{2}({\mathbb{R}}^{d})\rightarrow
L^{2}({\mathbb{R}}^{d})$ (see (\ref{opera}))  and its resolvent
$R_\lambda =(\mathcal{L} -\lambda )^{-1}$.  We will assume that $v
\geq 0$ is continuous, compactly supported and not identically
equal to zero. It is
well-known that the spectrum of $\mathcal{L}$ consists of the absolutely continuous part $%
(-\infty ,0]$ and at most a finite number of non-negative
eigenvalues:
\[
\sigma (\mathcal{L})=(-\infty ,0]\cup \{\lambda _{j}\},\text{ \ \
}0\leq j\leq N,\text{ \ \ }\lambda _{j} \geq 0.
\]
We enumerate the eigenvalues in the decreasing order. Thus, if
$\{\lambda _{j}\}\neq \emptyset $, then $\lambda _{0}=\max \lambda
_{j}$.  We assume that there is at least one positive eigenvalue.
The resolvent $R_\lambda: L^{2}({\mathbb{R}}^{d})\rightarrow
L^{2}({\mathbb{R}}^{d})$ is a meromorphic operator valued function
on $\mathbb{C}^{\prime }
 =\mathbb{C}\backslash (-\infty ,0]$.

 Denote the kernel of $R_\lambda$ by
$R_\lambda(x,y)$. If $v \equiv 0$ (in which case of course there
are no eigenvalues), the kernel depends on the difference $x-y$
and will intermittently use the notations $R^{0}_\lambda(x,y)$ and
$R^{0}_\lambda(x-y)$. The kernel $R^{0
}_\lambda(x)$ can be expressed through the Hankel function $%
H_{\nu }^{(1)}$:
\begin{equation}  \label{ha2}
R^{0}_\lambda(x)=
c_dk^{d-2}(k|x|)^{1-\frac{d}{2}}H_{\frac{d}{2}-1}^{(1)}(i\sqrt2k|x|),~~
k=\sqrt{\lambda},~~\mathrm{Re}k > 0.
\end{equation}

 We shall say that  $f \in C_{\rm \exp}( \mathbb{R}^d)$ (or simply
 $C_{\rm \exp}$)
if $f$ is continuous and
\[
||f||_{C_{\rm \exp}( \mathbb{R}^d)} = \sup_{x \in \mathbb{R}^d}(
|f(x)| e^{{|x|^2}}) < \infty.
\]
The space of bounded continuous functions on $\mathbb{R}^d$ will
be denoted  by $C( \mathbb{R}^d)$ or simply $C$. The following
simple lemma can be found in \cite{CKMV1}.
\begin{lemma} \label{rbetaW}
The operator $R_\lambda: C_{\exp}( \mathbb{R}^d) \rightarrow C(
\mathbb{R}^d)$ is meromorphic in $\lambda \in \mathbb{C}'$. Its
poles are of the first order and are located at eigenvalues of the
operator $\mathcal{L}$. For each $\varepsilon
> 0$ and some $\Lambda$, the operator is uniformly bounded
in $\lambda \in \mathbb{C}'$,
 $|{\rm arg} \lambda| \leq \pi -
\varepsilon$, $|\lambda| \geq \Lambda$. It is of order
$O(1/|\lambda|)$ as $\lambda \rightarrow \infty$, $|{\rm arg}
\lambda| \leq \pi - \varepsilon$. The eigenvalue $\lambda_0$ of
the operator $\mathcal{L}$ is simple and the corresponding
eigenfunction does not change sign.
\end{lemma}
From Lemma~\ref{rbetaW} it follows that the residue of $R_\lambda$
at $\lambda_0$ is the integral operator with the kernel $\psi(x)
\psi(y)$, where $\psi$ is the positive eigenfunction normalized by
the condition $||\psi||_{L^2( \mathbb{R}^d)} = 1$. The function
$\psi$ decays exponentially at infinity. More precisely, it
follows from (\ref{ha2}) that if we write $x$ as $(\theta, |x|)$
in polar coordinates, then there is a positive continuous function
$F$ such that
\begin{equation} \label{uuo}
\psi(x) \sim F(\theta) |x|^{\frac{1}{2}-\frac{d}{2}}
\exp(-\sqrt{2\lambda_0}|x|)~~~{\rm as}~~|x| \rightarrow \infty.
\end{equation}

 For a positive number $x$, we define the
curve $\Gamma(x)$ in the complex plane as follows:
\[
\Gamma(x) = \{\lambda: |{\rm Im} \lambda| = \sqrt{ 4x(x -{\rm Re}
\lambda)},~{\rm Re} \lambda \geq 0 \} \bigcup \{\lambda: |{\rm Im}
\lambda| = 2x(1 -{\rm Re} \lambda),~{\rm Re} \lambda \leq 0 \}.
\]
Thus $\Gamma(x)$ is a union of a piece of the parabola with the
vertex in $x$ that points in the direction of the negative real
axis and two rays tangent to the parabola at the points it
intersects the imaginary axis. The choice of the curve is somewhat
arbitrary, yet the following properties of $\Gamma(x)$ will be
important:

First, ${\rm Re} \lambda \leq x$ for $\lambda \in \Gamma(x)$.
Second, since the rays form a positive angle with the negative
real semi-axis, we have $|{\rm arg} \lambda| \leq \pi -
\varepsilon (x)$ for all $\lambda \in \Gamma(x)$ for some
$\varepsilon(x) > 0$. Third, since the rays are tangent to the
parabola, and the parabola is mapped into the line $\{ \lambda:
{\rm Re} \lambda = \sqrt{x} \}$ by the mapping $\lambda
\rightarrow \sqrt{\lambda}$, the image of the curve $\Gamma(x)$
under the same mapping lies in the half-plane $\{\lambda: {\rm Re}
\lambda \geq \sqrt{x} \}$.

 The integration along the vertical
lines in the complex plane and along contours $\Gamma(x)$, below,
is performed in the direction of the increasing complex part.

We'll need estimates on the solutions of the following parabolic
equation. Let
\begin{equation} \label{zeroind}
\partial_t \rho(t,x) = \frac{1}{2} \Delta \rho(t,x) +
 v(x) \rho(t,x),~~~~ \rho(0,x) = g(x) \in C_{\rm exp}.
\end{equation}
%\begin{equation} \label{zerorhs}
%\partial_t u(t,x) = \frac{1}{2} \Delta u(t,x) +
%\beta v(x) u(t,x) ,~~~~ u(0,x) =  g(x) \in C_{\rm exp}.
%\end{equation}
 We'll denote the Laplace transform of a function $f$ by $
\widetilde{f}$,
\[
\widetilde{f}(\lambda) =  \int_0^\infty \exp(-\lambda t) f(t) d t.
\]
%For a real number $x$, we denote the line $\{ \lambda: {\rm Re}
%(\lambda) = x \}$ in the complex plane  by $\Lambda_x$. The
%integration along $\Lambda_x$ is performed so that the imaginary
%part of the argument is increasing.
 Let $r$ be the distance between $\lambda_0$ and the rest of
the spectrum of the operator $\mathcal{L}$. In the arguments that
follow we'll use the symbol $A$ to denote constants that may
differ from line to line.

\begin{lemma} \label{tth}
For each $\varepsilon \in(0, r)$, the solution of (\ref{zeroind})
has the form
\begin{equation} \label{tthh}
\rho(t,x) = \exp(\lambda_0 t) \langle \psi,  g \rangle \psi(x) +
q(t,x),
\end{equation}
where
%\[
%C = -\int \psi_\beta(y) \int_0^\infty e^{-\lambda_0 t}
%\beta v(y)  f(t,y) d t,
%\]
\[
||q(t,\cdot)||_{C}  \leq A(\varepsilon)
\exp((\lambda_0-\varepsilon) t) ||g||_{C_{\rm exp}}.
\]
\end{lemma}
\proof After the Laplace transform, the equation becomes
\[
(\frac{1}{2} \Delta +  v) \widetilde{\rho}  - \lambda
\widetilde{\rho} = - g.
\]
Thus, the solution $\rho$ can be represented as
\begin{equation} \label{r11}
\rho(t,\cdot) = -\frac{1}{2 \pi i} \int_{{{\rm Re} \lambda =
\lambda_0 +1}} e^{\lambda t} R_\lambda  g d \lambda.
\end{equation}
The resolvent is meromorphic in the complex plane outside of the
interval $(-\infty, \lambda_0-r]$,  with the only (simple) pole
at~$\lambda_0$
 with the principal part of the Laurent expansion being
the integral operator with the kernel $\psi(x) \psi(y)/(\lambda_0
- \lambda)$.

By Lemma~\ref{rbetaW}, the norm of $R_\lambda$ does not exceed
$A/|\lambda|$ near infinity to the right of $\Gamma(\lambda_0
-\varepsilon)$. Therefore, the same integral as in (\ref{r11}) but
along the segment parallel to the real axis connecting a point
$\lambda_0 + 1 + i b$ with the contour $\Gamma(\lambda_0
-\varepsilon)$   tends to zero when $b \rightarrow \infty$.
Therefore, we can replace the contour of integration in
(\ref{r11}) by $\Gamma(\lambda_0 -\varepsilon)$. The residue gives
the main term, while
%\[
%-\exp(\lambda_0 t) \langle \psi_\beta, \widetilde{\beta v f}
%\rangle \psi_\beta(x).
%\]
% The
the integral over $\Gamma(\lambda_0 -\varepsilon)$ gives the
remainder term. \qed

\begin{lemma} \label{rhoo} Let $K \subset \mathbb{R}^d$ be a compact set. For
each $\varepsilon \in (0,r)$, the function $\rho_1(t,x,y)$
satisfies
\[
\rho_1(t,x,y) = \exp( \lambda_0 t) \psi(x) \psi(y)  + q_1 (t,x,y),
\]
where
\begin{equation} \label{equu}
\sup_{x \in K} |q_1(t,x,y)| \leq A(\varepsilon) \exp((\lambda_0-
\varepsilon) t - |y| \sqrt{2 (\lambda_0 - \varepsilon)})
\end{equation}
for $t \geq 1/2$. Moreover,
\begin{equation} \label{ort}
\int_{ \mathbb{R}^d} q_1(t,x,y)\psi(y) d y = 0
\end{equation}
for each $t > 0$, $x \in \mathbb{R}^d$, and
\begin{equation} \label{ep11}
\sup_{x \in \mathbb{R}^d} \int_{ \mathbb{R}^d} \rho_1(t,x,y) d y
\leq A \exp( \lambda_0 t).
\end{equation}
\end{lemma}
 \proof First, let us show that
\begin{equation} \label{iddt}
\langle \psi, \rho_1(t,\cdot,y) \rangle = \exp( {\lambda_0} t)
\psi(y)
\end{equation}
for each $t > 0$, $y \in \mathbb{R}^d$. Indeed,
\[
0 = \int_0^{t} \langle (\frac{\partial}{\partial s} +
\mathcal{L})(\exp(-\lambda_0 s) \psi), \rho_1 \rangle d s
=
\]
\[
\langle \exp(-\lambda_0 s) \psi, \rho_1 \rangle|_{s=0}^{t} +
\int_0^{t} \langle (\exp(-\lambda_0 s) \psi), (-
\frac{\partial}{\partial s} + \mathcal{L}) \rho_1 \rangle d s =
\]
\[
\langle \exp(-\lambda_0 t) \psi, \rho_1(t,\cdot,y) \rangle -
\langle  \psi, \rho_1 (0,\cdot,y) \rangle = \exp(-\lambda_0 t)
\langle
 \psi, \rho_1(t,\cdot,y) \rangle - \psi(y),
\]
which proves (\ref{iddt}). The relationship (\ref{iddt})
immediately implies (\ref{ort}).

Next, let $K'$ be a compact set that contains ${\rm supp}(v) \cup
K $ in its interior. In order to prove (\ref{equu}), consider
first the case when $y \in K'$. Apply (\ref{tthh}) with $t$
replaced by $t' = t - 1/2$ and $g = \rho_1(1/2,\cdot,y)$. In order
to calculate the main term of the asymptotics, we note that $
||g||_{C_{\rm exp}}$ is bounded uniformly in $y \in K'$ and
\[
\langle \psi, g \rangle = \exp(\frac{1}{2} {\lambda_0}) \psi(y),
\]
as follows from (\ref{iddt}). Therefore, (\ref{tthh}) implies that
\[
\rho_1(t,x,y) = \exp(\lambda_0 t)\psi(y) \psi(x) + \exp((\lambda_0
-\varepsilon) t ) q(t,x,y),
\]
where $||q(t,\cdot,y)||_C \leq A(K')$ for all $y \in K'$. We still
need to consider the case when $y \notin K'$.

Let $u(t,x,y) = \rho_1(t,x,y) - p_0(t,x,y)$, where $p_0$ is the
fundamental solution of the heat equation. Then $u$ satisfies the
non-homogeneous version of (\ref{zeroind}) with the right hand
side $f = - v(x) p_0(t,x,y)$ and $g \equiv 0$. Note that $f$ is a
smooth function since $y \notin K'$. Solving this equation for $u$
using the Laplace transform, as in the proof of Lemma~\ref{tth},
we obtain
\[
u(t,\cdot,y) = -\frac{1}{2 \pi i} \int_{{{\rm Re} \lambda =
\lambda_0 +1}} e^{\lambda t} R_\lambda  (- v
\widetilde{p}_0(\lambda, \cdot, y)) d \lambda
\]
\begin{equation} \label{iioo}
= - \frac{1}{2 \pi i} \int_{{{\rm Re} \lambda = \lambda_0 +1}}
e^{\lambda t} R_\lambda  ( v R^0_{\lambda}(\cdot,y)) d \lambda
\end{equation}
\[
=\exp(\lambda_0 t)\langle \psi,  v R^0_{\lambda_0}(\cdot,y)
\rangle \psi - \frac{1}{2 \pi i} \int_{ \Gamma( \lambda_0
-\varepsilon)} e^{\lambda t} R_\lambda  ( v
R^0_{\lambda}(\cdot,y)) d \lambda,
\]
where the first term on the right hand side is due to the residue
at $\lambda = \lambda_0$. The first term can be re-written as
\[
\exp(\lambda_0 t)\langle \psi,  v R^0_{\lambda_0}(\cdot,y) \rangle
\psi(x) =
\]
\[
\exp(\lambda_0 t)(R^0_{\lambda_0} ( v \psi)) (y) \psi(x) =
-\exp(\lambda_0 t)\psi(y) \psi(x).
\]
The last equality here follows from the fact that $\psi$ is an
eigenfunction with eigenvalue $\lambda_0$, that is
\[
(\frac{1}{2} \Delta  - \lambda_0) \psi = - v \psi.
\]
In order to estimate the second term on the right hand side of
(\ref{iioo}), we note that from~(\ref{ha2})  it follows that
\[
|R^{0}_\lambda(x,y)| \leq  A(l)
|\sqrt{\lambda}|^{\frac{d}{2}-\frac{3}{2}} |x - y
|^{\frac{1}{2}-\frac{d}{2}}  |\exp(-\sqrt{2 \lambda} |y-x|)|
\]
if $|\lambda|, |y -x| \geq l$. Thus
\[
 || v
R^0_{\lambda}(\cdot,y)||_{C_{\rm exp}} \leq A(\varepsilon) |y
|^{\frac{1}{2}-\frac{d}{2}}
|\sqrt{\lambda}|^{\frac{d}{2}-\frac{3}{2}} \exp(-\sqrt{2(\lambda_0
- \varepsilon)} |y|)
\]
for $y \notin K'$, $\lambda \in \Gamma(\lambda_0 - \varepsilon)$
due to the fact that ${\rm Re} \sqrt{\lambda} \geq \sqrt{\lambda_0
- \varepsilon}$ for $\lambda \in \Gamma(\lambda_0 - \varepsilon)$
and $|y - x| \geq l$ for $x \in {\rm supp}(v)$, $y \notin K'$.

Hence, using the estimate on the norm of $R_\lambda : C_{\rm exp}
\rightarrow C$ from Lemma~\ref{rbetaW}, we obtain
\[
|| R_\lambda  ( v R^0_{\lambda}(\cdot,y)) ||_C \leq A(\varepsilon)
|\sqrt{\lambda}|^{\frac{d}{2}-\frac{5}{2}} \exp(-\sqrt{2(\lambda_0
- \varepsilon)} |y|),~~~~\lambda \in \Gamma(\lambda_0 -
\varepsilon).
\]
Therefore, since ${\rm Re} \lambda \leq \lambda_0 -\varepsilon$
for $\lambda \in \Gamma(\lambda_0 -\varepsilon)$ and the factor
$e^{\lambda t}$ decays exponentially along $\Gamma(\lambda_0
-\varepsilon)$, the $C$-norm of the second term on the right hand
side of (\ref{iioo}) does not exceed $A(\varepsilon)
\exp((\lambda_0- \varepsilon) t - |y| \sqrt{2 (\lambda_0 -
\varepsilon)})$. The term $p_0(t,x,y)$ with $x \in K$, $y \notin
K'$, $t \geq 1/2$, is estimated by the same expression, possibly
with a different constant $A(\varepsilon )$. Indeed, if $t \geq
1/2$, then
\[
p_0(t,x,y) \leq A \exp(-|y-x|^2/2t) \leq A \exp((\lambda_0-
\varepsilon) t - |y-x| \sqrt{2 (\lambda_0 - \varepsilon)})
\]
since
\[
|y-x|^2/2t + (\lambda_0- \varepsilon) t - |y-x| \sqrt{2 (\lambda_0
- \varepsilon)} = (|y-x|/\sqrt{2t} - \sqrt{ (\lambda_0-
\varepsilon) t} ) ^2 \geq 0.
\]
This completes the proof of (\ref{equu}). In order to prove
(\ref{ep11}), we again write $u(t,x,y) = \rho_1(t,x,y) -
p_0(t,x,y)$. For fixed $x$, apply the Duhamel formula to the
equation
\[
\partial_t u(t,x,y) = \frac{1}{2} \Delta u(t,x,y) + v(y)\rho_1(t,x,y)
\]
with the initial data $u(0,x,\cdot) \equiv 0$.  Now (\ref{ep11})
follows since the $L^1$-norm of $v(y) \rho_1(t,x,y)$ is bounded by
$A \exp( \lambda_0 t)$ uniformly in $x$, as follows from
(\ref{equu}).
 \qed
 \\

We'll  need additional notations in order to describe the
asymptotics of $\rho_n$ with $n > 1$. Let
$\alpha^1_\varepsilon(t,y) = \psi(y)$ and
$\alpha^2_\varepsilon(t,y) = \exp(-\varepsilon t - |y| \sqrt{2
(\lambda_0 - \varepsilon)})$. Consider all possible sequences
$\sigma = (\sigma_1,...,\sigma_n)$ with $\sigma_i \in \{1,2\}$. By
$\Pi^n_\varepsilon (t,y_1,...,y_n)$ we denote the quantity
\[
\Pi^n_\varepsilon (t,y_1,...,y_n) = \sup_{\sigma \neq (1,...,1)}
\alpha^{\sigma_1}_\varepsilon(t,y_1)\cdot ... \cdot
\alpha^{\sigma_n}_\varepsilon(t,y_n).
\]
%Given functions $f_1(y_1)$, $f_2(y_1,y_2)$,...,
%$f_{n-1}(y_1,...,y_{n-1})$, we define
%\[
%F_n(f_1,...,f_{n-1}) (y_1,...,y_n) = \sum_{U \subset Y, U \neq
%\emptyset} f_{|U|} (U) f_{|Y \setminus U|}(Y\setminus U),
%\]
%where $Y = (y_1,...,y_n)$, $U$ is a proper non-empty subsequence
%of $Y$, and $|U|$ is the number of elements in this subsequence.

Let $P_t: C_{\rm \exp} \rightarrow C$ be the operator that maps
the initial function $g$ to the solution $\rho(t,\cdot)$ of
equation (\ref{zeroind}). Let $P^0_t g (x) = \exp(\lambda_0 t)
\langle \psi,  g \rangle \psi(x)$ and $P^1_t = P_t - P^0_t$.
Lemma~\ref{tth} states that
\[
||P^1_t|| \leq A(\varepsilon)
\exp((\lambda_0-\varepsilon) t).
\]
The particular form of $P^0_t$ then implies that
\begin{equation} \label{hjh}
||P_t|| \leq  ||P^0_t|| + ||P^1_t|| \leq A'  \exp(\lambda_0 t).
\end{equation}
For $g \in C_{\rm exp}$ and $n \geq 2$,  we denote
\[
I_n(g) := R_{ n \lambda_0} g  = \int_0^\infty \exp(-n \lambda_0 s)
P_s g d s \in C.
\]
Note that
\[
\int_0^t  \exp(n \lambda_0 s) P_{t-s} g d s = \exp(n \lambda_0 t)
\int_0^t  \exp(-n \lambda_0 s) P_s g d s
\]
\begin{equation} \label{amm}
= \exp(n \lambda_0 t)(I_n(g) + O(\exp(-(n-1) \lambda_0 t))
)~~~~{\rm as}~~t \rightarrow \infty.
\end{equation}
The functions $f_1, f_2,...$  are defined inductively: $f_1 =
\psi$ and
\[
f_n =   \sum_{k=1}^{n-1} \frac{n!}{k! (n-k)!} I_n(v f_k
f_{n-k}),~~~n \geq 2.
\]

\begin{lemma} \label{rnn}
Let $K \subset \mathbb{R}^d$ be a compact set. For each
$\varepsilon \in (0,r)$, the function $\rho_n$ satisfies
\begin{equation} \label{yuy2}
\rho_n(t,x,y_1,...,y_n) =  \exp( n \lambda_0 t) f_n(x) \psi(y_1)
\cdot ... \cdot \psi(y_n) + q_n(t,x,y_1,...,y_n),
\end{equation}
where
\begin{equation} \label{ggp}
\sup_{x \in K} |q_n(t,x,y_1,...,y_n)| \leq A_n(\varepsilon) \exp(n
\lambda_0 t ) \Pi^n_\varepsilon (t,y_1,...,y_n)
\end{equation}
for $t \geq 1/2$.
\end{lemma}
\proof For $n =1$, the relation (\ref{yuy2}) coincides with the
statement of Lemma~\ref{rhoo}. Let us assume that (\ref{yuy2})
holds for all natural numbers up to and including $n-1$. A generic
subsequence $U \subset Y = (y_1,...,y_n)$ will be written as $U =
(z_1,...,z_{|U|})$ and its complement as $Y \setminus U =
(\overline{z}_1,...,\overline{z}_{n - |U|})$. By the Duhamel
principle applied to the equation for $\rho_n$, we obtain
\[
\rho_n(t,\cdot,y_1,...,y_n) = \int_0^t P_{t-s} ( v \sum_{U \subset
Y, U \neq \emptyset} \rho_{|U|} (s,\cdot,z_1,...,z_{|U|})
\rho_{n-|U|}(s,\cdot,\overline{z}_1,...,\overline{z}_{n - |U|})) d
s
\]
\[
= \int_0^t P_{t-s} ( v  \sum_{U \subset Y, U \neq \emptyset} \exp(
|U| \lambda_0 s) f_{|U|} (\cdot) \psi(z_1)\cdot...\cdot
\psi(z_{|U|})
\]
\begin{equation} \label{middxx}
\times  \exp( (n-|U|) \lambda_0 s) f_{n - |U|} (\cdot)
\psi(\overline{z}_1)\cdot...\cdot \psi(\overline{z}_{n-|U|})) d s
\end{equation}
\[
+  2 \int_0^t P_{t-s} ( v  \sum_{U \subset Y, U \neq \emptyset}
\exp( |U| \lambda_0 s) f_{|U|} (\cdot) \psi(z_1)\cdot...\cdot
\psi(z_{|U|}) q_{n-|U|}(s,\cdot,\overline{z}_1,...,\overline{z}_{n
- |U|})) d s
\]
\[
+ \int_0^t P_{t-s} ( v  \sum_{U \subset Y, U \neq \emptyset}
q_{|U|}(s,\cdot,{z}_1,...,{z}_{|U|})
q_{n-|U|}(s,\cdot,\overline{z}_1,...,\overline{z}_{n - |U|})) d s.
\]
The second and third integrals on the right hand side of
(\ref{middxx}) contribute only to the remainder term. Indeed,
consider the contribution to the second integral from the term
with a given $U$:
\[
\int_0^t P_{t-s} ( v  \exp( |U| \lambda_0 s) f_{|U|} (\cdot)
\psi(z_1)\cdot...\cdot \psi(z_{|U|})
q_{n-|U|}(s,\cdot,\overline{z}_1,...,\overline{z}_{n - |U|})) d s
\]
\[
\leq A \psi(z_1)... \psi(z_{|U|}) \int_0^t P_{t-s} ( v \exp( |U|
\lambda_0 s) f_{|U|} (\cdot)
\]
\[
\times
 \exp((n - |U|)
\lambda_0 s ) \Pi^{n-|U|}_\varepsilon
(s,\overline{z}_1,...,\overline{z}_{n - |U|})) d s
\]
\[
\leq  A \psi(z_1)\cdot...\cdot \psi(z_{|U|}) \int_0^t
\exp(\lambda_0(t-s)) \exp(n \lambda_0 s) \Pi^{n-|U|}_\varepsilon
(s,\overline{z}_1,...,\overline{z}_{n - |U|})) d s
\]
\[
\leq  A \exp(n \lambda_0 t) \psi(z_1)... \psi(z_{|U|})
\Pi^{n-|U|}_\varepsilon (t,\overline{z}_1,...,\overline{z}_{n -
|U|}) \leq  A\exp(n \lambda_0 t ) \Pi^n_\varepsilon
(t,y_1,...,y_n),
\]
where the first inequality follows from the inductive assumption
and the second one from~(\ref{hjh}). The third integral on the
right hand side of (\ref{middxx}) is estimated similarly. It
remains to consider the first integral. It is equal to
\[
\psi(y_1)\cdot...\cdot \psi(y_n)
 \int_0^t   \exp (n \lambda_0 s) P_{t-s} (
v    \sum_{U \subset Y, U \neq \emptyset} f_{|U|} f_{n-|U|} ) d s
\]
\[
= \psi(y_1) \cdot ... \cdot \psi(y_n)   \exp( n \lambda_0 t)
\left( f_n(\cdot) +  O(\exp(-(n-1) \lambda_0 t))\right),
\]
where the last equality follows from (\ref{amm}). Thus we obtain
the main term from the right hand side of (\ref{yuy2}) plus the
correction
\[
\psi(y_1) \cdot ... \cdot \psi(y_n)   \exp( n \lambda_0 t)
O(\exp(-(n-1) \lambda_0 t))
\]
for which the estimate (\ref{ggp}) holds since $ \psi(y_1)
\exp(-\lambda_0 t) \leq  \alpha^2_\varepsilon(t,y_1)$ due to
(\ref{uuo}). \qed
\\

\section{Growth of the total number of particles} \label{growth}

We denote the probability space on which the branching process is
defined by $(\Omega, \mathcal{F}, \mathrm{P})$. Let
$\mathcal{F}_t$, $t \geq 0$,  be the filtration generated by the
process. We'll write $L^2$ for $L^2(\Omega, \mathcal{F},
\mathrm{P})$. Let $N^x_t$ be the number of particles in $
\mathbb{R}^d$ at time $t$, assuming that at $t = 0$ there was a
single particle located at $x$.

In this section we prove the basic result on the convergence of
${N^x_t}/e^{\lambda_0 t}$ in $L^2$. The almost sure convergence
and the asymptotics of the number of particles in a (possibly
time-dependent) region of space will be considered in the
following sections.
\begin{theorem} \label{ftc}
There is a random variable $\xi^x$ such that
\begin{equation} \label{mli}
\frac{N^x_t}{e^{\lambda_0 t}}  \rightarrow \xi^x~~{\it as}~~ t
\rightarrow \infty,
\end{equation}
where the convergence takes place in $L^2$.
\end{theorem}
\proof Observe that for $0 < s \leq t$,
\begin{equation} \label{product}
\mathrm{E} (N^x_s N^x_t) =  \int_{ \mathbb{R}^d}   \int_{
\mathbb{R}^d}  \int_{ \mathbb{R}^d} \rho_2 (s, x, y_1, z)
\rho_1(t-s, z, y_2) d z d y_1 d y_2 + \int_{ \mathbb{R}^d}
\rho_1(t, x, y_2) dy_2.
\end{equation}
Indeed, fix $y_1, y_2 \in \mathbb{R}^d$. Then the probability that
there is a particle in an infinitesimal neighborhood of $y_1$ at
time $s$, while a different particle present at time $s$ gives
rise to a particle in an infinitesimal neighborhood of $y_2$ at
time $t$ is equal to
\[
\left(\int_{ \mathbb{R}^d} \rho_2 (s, x, y_1, z) \rho_1(t-s, z,
y_2) d z \right) d y_1 d y_2.
\]
The probability that a particle in an infinitesimal neighborhood
of $y_1$ at time $s$ gives rise to a particle in an infinitesimal
neighborhood of $y_2$ at time $t$ is equal to
\[
\rho_1(s,x,y_1) \rho_1(t-s, y_1, y_2)d y_1 d y_2.
\]
After adding the contributions from the two events and integrating
in $y_1$ and $y_2$, we obtain~(\ref{product}).

Combining (\ref{product}) with (\ref{yuy2}) and using that $f_1 =
\psi$, we see that
\[
\mathrm{E} (N^x_s N^x_t) =  \int_{ \mathbb{R}^d}   \int_{
\mathbb{R}^d}  \int_{ \mathbb{R}^d}  [ e^{2 \lambda_0 s} f_2(x)
\psi(y_1) \psi(z)  e^{\lambda_0(t-s)} \psi(z) \psi(y_2) +
\]
\[
e^{2 \lambda_0 s} f_2(x) \psi(y_1) \psi(z)  q_1(t-s,z, y_2) +
q_2(s,x,y_1,z) \rho_1(t-s,z,y_2) ]  d z d y_1 d y_2 +
\]
\[
 \int_{ \mathbb{R}^d}
\rho_1(t, x, y_2) dy_2 =: I^1_{s,t}(x) + I^2_{s,t}(x) +
I^3_{s,t}(x) + I^4_{s,t}(x).
\]
Note that
\[
I^1_{s,t} = e^{\lambda_0(s+t)} f_2(x) (\int_{ \mathbb{R}^d } \psi
)^2
\]
since $\int_{ \mathbb{R}^d } \psi^2(z) d z = 1$.  Also observe
that $ I^2_{s,t}(x) = 0$ since $\int_{ \mathbb{R}^d}  \psi(z)
q_1(t-s,z, y_2) d z = 0$ by (\ref{ort}). Finally,
\[
\sup_{x \in K} |I^3_{s,t}(x) + I^4_{s,t}(x)| \leq A
e^{\lambda_0(t + s) - \varepsilon s}
\]
by Lemma~\ref{rnn} and (\ref{ep11}). Therefore,
\[
\sup_{x \in K} |\mathrm{E} (N^x_s N^x_t) - e^{\lambda_0(s+t)}
f_2(x) (\int_{ \mathbb{R}^d } \psi )^2 |  \leq A e^{\lambda_0(t +
s) - \varepsilon s}.
\]
Thus we have
\[
\mathrm{E} \left(\frac{N^x_s}{e^{\lambda_0 s}} -
\frac{N^x_t}{e^{\lambda_0 t}}\right)^2 = \frac{\mathrm{E}
(N^x_s)^2}{e^{2 \lambda_0 s}} + \frac{\mathrm{E} (N^x_t)^2}{e^{2
\lambda_0 t}} - 2 \frac{\mathrm{E} (N^x_s N^x_t)}{e^{\lambda_0
(s+t)}} \leq A e^{-\varepsilon s}
\]
for $x \in K$. This shows that  $ {N^x_t}/{e^{\lambda_0 t}}$ is a
Cauchy family of random variables as $t \rightarrow \infty$, and
we have convergence in $L^2$. \qed
\\
\\
{\bf Remark.} It is possible to show (see \cite{K}) that all the
moments of the variables $ {N^x_t}/{e^{\lambda_0 t}}$ converge to
those of $\xi^x$. The moments of the limiting distribution are
\[
{\mathrm{E}}(\xi^{x})^n = \left( \int_{ \mathbb{R}^d} \psi(y) d y
\right)^n f_n(x).
\]
They  were shown to determine the distribution of $\xi^x$
uniquely.
\\
\\
\noindent
 {\bf Remark.} In dimensions $d \geq 3$ the limiting
random variable $\xi^x$ is equal to zero with positive
probability. Indeed, since the diffusion is transient, there is a
positive probability that the original particle wanders off to
infinity without branching. Let $ B^x$ be the event that the
number of particles stays bounded (and therefore tends to a finite
limit as $t \rightarrow \infty$) and $ E^x  =  \{\xi^x > 0 \}$ be
the event that the number of particles grows exponentially. It is
possible to show (see \cite{K}, for example) that $ \mathrm{P}(B^x
\cup E^x) = 1$ for each $x$.

\section{Growth of the number of particles in a domain} \label{growth22}

Let $n^{x}_t(U)$ denote the number of particles in a domain $U
\subseteq \mathbb{R}^d$, assuming that at $t = 0$ there was a
single particle located at $x$. Let
\[
\alpha(U) = \frac{\int_U \psi(y) d y }{ \int_{ \mathbb{R}^d}
\psi(y) d y}~.
\]
The asymptotics of the number of particles in $U$ is given by the
following theorem.
\begin{theorem} \label{domaint} For each measurable $U \subseteq
\mathbb{R}^d$, we have
\begin{equation} \label{domainconv}
 \frac{n^{x}_t(U) }{e^{\lambda_0 t}}  \rightarrow
\alpha(U) \xi^x ~~{\it as}~~ t \rightarrow \infty,
\end{equation}
where the convergence takes place in $L^2$.
\end{theorem}
\proof  By Theorem~\ref{ftc} it is sufficient to prove that
\[
 \frac{n^{x}_t(U) - \alpha(U) N^x_t }{e^{\lambda_0 t}}  \rightarrow
0~~{\rm as}~~ t \rightarrow \infty.
\]
Observe that
\begin{equation} \label{oo1}
\mathrm{E} (N^x_t)^2 =  \int_{ \mathbb{R}^d}   \int_{\mathbb{R}^d}
 \rho_2 (t, x, y_1, y_2)  dy_1 d y_2 + \int_{
\mathbb{R}^d} \rho_1(t, x, y_1) dy_1,
\end{equation}
\begin{equation} \label{oo2}
\mathrm{E} (n^x_t(U) N^x_t) =  \int_{ \mathbb{R}^d}   \int_{U}
 \rho_2 (t, x, y_1, y_2)  dy_1 d y_2 + \int_{
U} \rho_1(t, x, y_1) dy_1,
\end{equation}
\begin{equation} \label{oo3}
\mathrm{E} (n^x_t(U))^2 =  \int_{ U}   \int_{U}
 \rho_2 (t, x, y_1, y_2)  dy_1 d y_2 + \int_{
U} \rho_1(t, x, y_1) dy_1.
\end{equation}
Upon expanding $(n^{x}_t(U) - \alpha(U) N^x_t )^2 = (n^{x}_t(U))^2
- 2\alpha(U) n^{x}_t(U) N^x_t + (\alpha(U)N^x_t)^2$, using
Lemma~\ref{rnn} for the asymptotics of $\rho_1$ and $\rho_2$ and
collecting all the lower order terms in the remainder $R(t,x)$, we
obtain
\[
\mathrm{E} \left( \frac{n^{x}_t(U) - \alpha(U) N^x_t
}{e^{\lambda_0 t}} \right)^2  =
\]
\[
f_2(x) ( \int_{ U}   \int_{U} \psi(y_1) \psi(y_2) d y_1 dy_2 - 2
\alpha(U) \int_{ \mathbb{R}^d}   \int_{U} \psi(y_1) \psi(y_2) d
y_1 dy_2 +
\]
\[
(\alpha(U))^2  \int_{ \mathbb{R}^d} \int_{\mathbb{R}^d} \psi(y_1)
\psi(y_2) d y_1 dy_2 ) + R(t,x) =
\]
\[
f_2(x) \left(\int_U \psi - \alpha(U) \int_{ \mathbb{R}^d} \psi
\right)^2 + R(t,x) = R(t,x),
\]
where $\sup_{x \in K} R(t,x) \leq A e^{-\varepsilon t}$. This
shows that $ (n^{x}_t(U) - \alpha(U) N^x_t )/{e^{\lambda_0 t}}$ is
a Cauchy family, thus completing the proof. \qed
\\

Now let us examine the number of particles in the vicinity of a
point that is at a linear in $t$ distance from the origin. Let $b
= \sqrt{\lambda_0/2}$ ($b$ will be seen to be the rate of growth
of the region where the particles can be found with probability
that tends to one). Let $v \in  \mathbb{R}^d$ be a vector with $0
<  |v| < b$ and $U_{tv} = U + t v $ the domain obtained from $U$
by translation by the vector $t v$. Let
\[
g(t) = g(U, t, v) =  \frac{ e^{\lambda_0 t}  \int_{U_{tv}} \psi(y)
d y}{ \int_{ \mathbb{R}^d} \psi(y) d y}.
\]
Note that the asymptotics of $g(t)$ can be obtained from
(\ref{uuo}). In particular,
\begin{equation} \label{ggee}
A_1 t^{\frac{1-d}{2}} e^{(\lambda_0 - \sqrt{2 \lambda_0} |v|) t}
\leq g(t) \leq A_2 t^{\frac{1-d}{2}} e^{(\lambda_0 - \sqrt{2
\lambda_0} |v|) t}
\end{equation}
if $U$ is bounded.

%For $c \in (0,b)$ and $\theta \in \mathbb{R}^d$, $|\theta| = 1$,
%we shall study the number of particles in $B_{ c t \theta }$, the
%ball of radius one centered at $ c t \theta$.
%
%Let
%\[
%\beta(\theta) = \frac{  F(\theta) \int_{B_0} e^{y_1 \sqrt{2
%\lambda_0} } d y}{ \int_{ \mathbb{R}^d} \psi(y) d y},
%\]
%where $F$ is the same as in~(\ref{uuo}), and define
%\[
%g(t) = (c t)^{\frac{1-d}{2}} e^{(\lambda_0 - \sqrt{2 \lambda_0} c)
%t}.
%\]
\begin{theorem} \label{domaint2} Let $U$ be a bounded domain. For each $v \in \mathbb{R}^d$ such that
$|v| < b$, we have
\begin{equation} \label{domainconv2}
 \frac{n^{x}_t(U_{tv}) }{g(t)}  \rightarrow \xi^x ~~{\it as}~~ t \rightarrow \infty,
\end{equation}
where the convergence takes place in $L^2$.
\end{theorem}
\proof  From (\ref{oo1})-(\ref{oo3}) we obtain
\begin{equation} \label{rha}
\mathrm{E} \left( \frac{n^{x}_t(U_{tv}) }{g(t)} -
\frac{N^x_t}{e^{\lambda_0 t}} \right)^2 = f_2(x) \left(
\frac{e^{\lambda_0 t}}{g(t)}  \int_{ U_{tv}} \psi - \int_{
\mathbb{R}^d} \psi \right)^2 + R(t,x),
\end{equation}
where
\[
R(t,x) = ( \frac{1}{g^2(t)} \int_{U_{tv}}  \int_{ U_{tv}}
q_2(t,x,y_1,y_2) dy_1 dy_2 - \frac{ 2}{g(t) e^{\lambda_0 t}}
\int_{ U_{tv}} \int_{ \mathbb{R}^d} q_2(t,x,y_1,y_2) dy_1 dy_2 +
\]
\[
\frac{1}{e^{2 \lambda_0 t}} \int_{ \mathbb{R}^d} \int_{
\mathbb{R}^d} q_2(t,x,y_1,y_2) dy_1 dy_2 + \frac{1}{g^2(t)} \int_{
U_{tv}} \rho_1(t,x,y_1) dy_1
-
\]
\[
\frac{ 2}{g(t) e^{\lambda_0 t}} \int_{ U_{tv}} \rho_1(t,x,y_1)
dy_1 +  \frac{1}{e^{2 \lambda_0 t}} \int_{ \mathbb{R}^d}
\rho_1(t,x,y_1) dy_1).
\]
The first term on the right hand side of (\ref{rha}) is equal to
zero as follows from the definition of $g(t)$. Applying
Lemma~\ref{rnn} to estimate the integrand in each of the terms in
$R(t,x)$, it is easy to see that $R(t,x)$ tends to zero
exponentially fast as $t \rightarrow \infty$. Indeed, let us
consider the first term. By (\ref{ggee}) and (\ref{ggp}), with an
arbitrary $\varepsilon \in (0,r)$, we have the estimate
\[
|\frac{1}{g^2(t)} \int_{U_{tv}}  \int_{U_{tv}} q_2(t,x,y_1,y_2)
dy_1 dy_2| \leq
\]
\begin{equation} \label{rsti}
A(\varepsilon)\frac{  e^{2 \lambda_0 t - \varepsilon t - t |v|
\sqrt{2(\lambda_0 - \varepsilon)}}(e^{-\varepsilon t - t |v|
\sqrt{2 (\lambda_0 - \varepsilon)}} +  t^{\frac{1-d}{2}} e^{- t
|v| \sqrt{2 \lambda_0}} ) }{t ^{1-d} e^{2(\lambda_0 - \sqrt{2
\lambda_0} |v|) t} } \leq
\end{equation}
\[
A(\varepsilon) ( t^{d-1} e^{2t\left(|v|\sqrt{2 \lambda_0}  -
\varepsilon - |v|\sqrt{2(\lambda_0 - \varepsilon)}\right)} +
t^{\frac{d-1}{2}} e^{t\left(|v|\sqrt{2 \lambda_0}  - \varepsilon -
|v| \sqrt{2(\lambda_0 - \varepsilon)}\right)}).
\]
Note that $|v|\sqrt{2 \lambda_0}  - \varepsilon - |v|
\sqrt{2(\lambda_0 - \varepsilon)} < 0$ since $|v| \in (0,
\sqrt{\lambda_0/2})$. Therefore, the right hand side of
(\ref{rsti}) tends to zero exponentially fast as $t \rightarrow
\infty$. The other terms in the expression for $R(t,x)$ can be
dealt with in the same fashion.
  \qed
\\
\\
{\bf Remark.} With the help of Lemma~\ref{rnn} it is possible to
show that we have the convergence of all the moments in
(\ref{domainconv}) and (\ref{domainconv2}).

\section{The almost sure convergence} \label{asco}

\begin{theorem} \label{alms}
The convergence in  (\ref{mli})  takes place almost surely. If $U$
is a domain with a smooth boundary, then the convergence in
(\ref{domainconv}) and (\ref{domainconv2}) takes place almost
surely.
\end{theorem}
\proof We will only prove the almost sure convergence in
(\ref{domainconv2}) since the other statements can be proved
similarly. Fix an arbitrary $\delta > 0$ and $K > 0$.
%We need to demonstrate
%that there is $T = T(\delta)$ such that
%\[
%\mathrm{P}( | \frac{n^{x}_t(U_{tv}) }{g(t)} - \xi^x| > \delta
%~~~{\rm for}~~{\rm some}~~t \geq T ) \leq \delta.
%\]
By the Borel-Cantelli lemma, it is sufficient to demonstrate that
there is an increasing sequence $t_n \rightarrow \infty$ such that
\[
\sum_{n  = 1}^\infty \mathrm{P}(\sup_{t \in [t_n, t_{n+1}]} |
\frac{n^{x}_t(U_{tv}) }{g(t)} - \xi^x| > \delta,~\xi^x < K ) <
\infty.
\]
From the proof of Theorem~\ref{domaint2} it follows that $
{n^{x}_t(U_{tv}) }/{g(t)}$ converges to $\xi^x$ in $L^2$
exponentially fast. Let $\gamma > 0$ be such that
\begin{equation} \label{ecooi}
\mathrm{E}(\frac{n^{x}_t(U_{tv}) }{g(t)} - \xi^x)^2 \leq  A e^{-
\gamma t}
\end{equation}
for some constant $A$. We take $t_n = 2 \ln n/\gamma$, $n \geq 1$.
By the Chebyshev inequality,
\[
\mathrm{P}(| \frac{n^{x}_{t_n}(U_{t_nv}) }{g(t_n)} - \xi^x| >
\frac{\delta}{5} ) \leq \frac{25 A  e^{- \gamma t_n}}{\delta^2} =
\frac{25 A }{\delta^2 n^2},
\]
and therefore
\begin{equation} \label{jj11}
\sum_{n  = 1}^\infty \mathrm{P}(| \frac{n^{x}_{t_n}(U_{t_nv})
}{g(t_n)} - \xi^x| > \frac{\delta}{5} ) < \infty.
\end{equation}
It remains to show that
\[
\sum_{n  = 1}^\infty \mathrm{P}(\sup_{t \in [t_n, t_{n+1}]} |
\frac{n^{x}_t(U_{tv}) }{g(t)} - \frac{n^{x}_{t_n}(U_{t_nv})
}{g(t_n)}| > \frac{4\delta}{5},~\xi^x < K ) < \infty,
\]
which is equivalent to the following two inequalities holding at
the same time
\begin{equation} \label{kkj1x} \sum_{n  =
1}^\infty \mathrm{P}(\sup_{t \in [t_n, t_{n+1}]}
\frac{n^{x}_t(U_{tv}) }{g(t)} - \frac{n^{x}_{t_n}(U_{t_nv})
}{g(t_n)} > \frac{4\delta}{5},~\xi^x < K ) < \infty,
\end{equation}
\begin{equation} \label{kkj1y} \sum_{n  =
1}^\infty \mathrm{P}(\sup_{t \in [t_n, t_{n+1}]} (-
\frac{n^{x}_t(U_{tv}) }{g(t)} + \frac{n^{x}_{t_n}(U_{t_nv})
}{g(t_n)}) > \frac{4\delta}{5},~\xi^x < K ) < \infty.
\end{equation}
We will only prove (\ref{kkj1x}) since (\ref{kkj1y}) can be proved
similarly. Let us show that the term ${n^{x}_{t_n}(U_{t_nv})
}/{g(t_n)}$ in (\ref{kkj1x}) can be replaced by a more convenient
expression. For $r > 0$, let $U^r = \{x \in \mathbb{R}^d: {\rm
dist} (x, U) < r \}$ be the $r$-neighborhood of $U$. Let $g^r(t) =
g(U^r, t, v)$. Since $U$ is a smooth domain, from the definition
of $g$ it follows that $r$ can be chosen to be sufficiently small
so that
\begin{equation} \label{egp}
\sup_{t \in [t_n, t_{n+1}]} | \frac{{g^r}(t_n) - g(t) }{g(t)}| <
\frac{\delta}{5}/(K+ \frac{\delta}{5})
\end{equation}
for all sufficiently large $n$. Moreover, since $r >0$ and
$t_{n+1} - t_n \rightarrow \infty$, we have
\begin{equation} \label{inside1}
\bigcup_{t \in [t_n, t_{n+1}]} U_{tv} \subset U^{\frac{r}{2}}_{t_n
v}
%~~~~{\rm with}~~~~{\rm dist} (\bigcup_{t \in [t_n, t_{n+1}]}
%U_{tv},~ \mathbb{R}^d \setminus U^r_{t_n v}) \geq \frac{r}{2}
\end{equation}
for all sufficiently large $n$. As in (\ref{jj11}), we have
\begin{equation} \label{jj112}
\sum_{n  = 1}^\infty \mathrm{P}(| \frac{n^{x}_{t_n}(U^r_{t_nv})
}{g^r(t_n)} - \xi^x| > \frac{\delta}{5} ) < \infty.
\end{equation}
Now,
\[
\sum_{n  = 1}^\infty \mathrm{P}(\sup_{t \in [t_n, t_{n+1}]}
 | \frac{n^{x}_{t_n}(U^r_{t_nv})}{g(t)} -
\frac{n^{x}_{t_n}(U_{t_nv})}{g(t_n)}| > \frac{3\delta}{5},~\xi^x <
K ) \leq
\]
\[
\sum_{n  = 1}^\infty \mathrm{P}(\sup_{t \in [t_n, t_{n+1}]}
 | \frac{n^{x}_{t_n}(U^r_{t_nv})}{g^r(t_n)} -
\frac{n^{x}_{t_n}(U_{t_nv})}{g(t_n)}| > \frac{2\delta}{5},~\xi^x <
K ) +
\]
\[
\sum_{n  = 1}^\infty \mathrm{P}(\sup_{t \in [t_n, t_{n+1}]}
 | \frac{n^{x}_{t_n}(U^r_{t_nv})}{g(t)} -
\frac{n^{x}_{t_n}(U^r_{t_nv})}{g^r(t_n)}| >
\frac{\delta}{5},~\xi^x < K ).
\]
The first series on the right hand side is finite by (\ref{jj11})
and (\ref{jj112}). The second series is estimated from above by
\[
\sharp \{n: \sup_{t \in [t_n, t_{n+1}]} (K+ \frac{\delta}{5})|
\frac{g^r(t_n) - g(t) }{g(t)}|  > \frac{\delta}{5} \} + \sum_{n  =
1}^\infty \mathrm{P}(| \frac{n^{x}_{t_n}(U^r_{t_nv}) }{g^r(t_n)} -
\xi^x| > \frac{\delta}{5} ).
\]
The first term of this expression is finite as follows from
(\ref{egp}), while the second one is finite by (\ref{jj112}).
Therefore, (\ref{kkj1x}) will follow if we demonstrate that
\begin{equation} \label{sty1}
\sum_{n  = 1}^\infty \mathrm{P}(\sup_{t \in [t_n, t_{n+1}]} \frac{
{n^{x}_t(U_{tv}) } - {n^{x}_{t_n}(U^r_{t_nv}) }}{g(t)}
> \frac{\delta}{5},~\xi^x < K ) < \infty.
\end{equation}
Similarly to (\ref{jj11}), we have
\[
\sum_{n  = 1}^\infty \mathrm{P}(| \frac{N^{x}_{t_n} }{e^{\lambda_0
t_n}} - \xi^x| > c) < \infty
\]
for each constant $c >0$. Combining this with (\ref{ggee}) and
(\ref{inside1}), we see that (\ref{sty1}) will follow if we show
that for each $\alpha, R
> 0$ we have
\begin{equation} \label{sty1x}
\sum_{n  = 1}^\infty \mathrm{P}(\sup_{t \in [t_n, t_{n+1}]}
{n^{x}_t(U^{\frac{r}{2}}_{t_nv}) } - {n^{x}_{t_n}(U^r_{t_nv}) }
> n^\alpha,~ N^{x}_{t_n} < R e^{\lambda_0 t_n}  ) < \infty.
\end{equation}
Roughly speaking, we need to show that  the number of particles
that visit a smaller region $ U^{\frac{r}{2}}_{t_nv}$ over a short
interval of time $[t_n, t_{n+1}]$ can't significantly exceed the
number of particles that are in the larger region $ U^r_{t_nv}$ at
the time $t_n$. Let us fix $n$ and study the $n$-th term in the
series (\ref{sty1x}). After conditioning on $\mathcal{F}_{t_n}$,
the questions becomes the following: Suppose we have $m \leq R
e^{\lambda_0 t_n}$ particles at time zero located at
$x_1,...,x_m$. Let $y_1,...,y_M$ be their descendants at time $t =
t_{n+1} - t_n$. Each of these has a starting point in the set
$\{x_1,...,x_m\}$. We are interested in the probability that at
least $n^\alpha$ descendants of the original particles were at a
distance $r/2$ from their respective starting points prior to the
time $t = t_{n+1} - t_n$.

The expected number of descendants of a single particle that cover
distance $r/2$ (from the initial position of the particle) in time
$t$ decays faster than $e^{-\beta/t}$ as $t \rightarrow 0$ for
some $\beta > 0$. Therefore, by the Chebyshev inequality, the
$n$-th term in (\ref{sty1x}) is estimated from above by $ { R
e^{\lambda_0 t_n} e^{-\beta/(t_{n+1} - t_n)} } n^{-\alpha}$, which
decays exponentially in $n$, as follows from the definition of
$t_n$. Therefore the series (\ref{sty1x}) converges, which
completes the proof.
 \qed
\\

\section{Limiting shape of the region occupied by particles}
\label{shape}

Let $B(r)$ denote the ball of radius $r$ centered at the origin.
Recall that $b = \sqrt{\lambda_0/2}$.
\begin{theorem}
For each $\delta > 0$, there exists a random variable $T =
T(\delta)$ ($T < \infty$ almost surely) with the following
properties:

(a) There are no particles outside $B((b+\delta) t)$ for $t \geq
T$.

(b) On the event $\xi^x > 0$ the union of the unit neighborhoods
of the particles cover $B((b-\delta) t)$  for all $t \geq T$.
\end{theorem}
\noindent
 {\it Sketch of the proof.} We'll only verify the statement for a
sequence of times $t_n = c \ln n$ for a certain $c$. Namely, we'll
show that there is a random variable $N$ and a constant $c > 0$
such that:

$(a')$ There are no particles outside of $B((b+\delta) t_n)$ for
$n \geq N$.

$(b')$ On the event $\xi^x > 0$  the union of the unit
neighborhoods of the particles cover $B((b-\delta) t_n)$ for $n
\geq N$.

The transition from the sequence of times to the continuous time
can be then accomplished similarly to the way it was done in the
proof of Theorem~\ref{alms}.

By the Chebyshev inequality,
\[
\mathrm{P}( {n^{x}_{t_n}(B((b+\delta) t_n)) } \geq 1) \leq
\mathrm{E} {n^{x}_{t_n}(B((b+\delta) t_n)) } = \int_{ B((b+\delta)
t_n)} \rho_1(t_n, x, y) d y \leq
\]
\[
A_1 \int_{ B((b+\delta) t_n)} (\exp(\lambda_0 t_n {-|y|\sqrt{2
\lambda_0} }) + \exp((\lambda_0- \varepsilon) t_n - |y| \sqrt{2
(\lambda_0 - \varepsilon)}) ) d y \leq
\]
\[
A_2 \left( \exp( t_n (\lambda_0 - \sqrt{2 \lambda_0}( b +
\delta))) + \exp( t_n (\lambda_0 - \varepsilon - (b + \delta)
\sqrt{2(\lambda_0 - \varepsilon)}) ) \right),
\]
where the second inequality is due to Lemma~\ref{rhoo} and
(\ref{uuo}). By choosing a sufficiently small $\varepsilon > 0$,
we can make the right hand side of the last formula smaller than
$A_2 \exp(-t_n \delta \sqrt{\lambda_0})$. Therefore,
\[
\sum_{n  = 1}^\infty \mathrm{P}( {n^{x}_{t_n}(B((b+\delta) t_n)) }
\geq 1)  \leq A_2\sum_{n  = 1}^\infty e^{-t_n \delta
\sqrt{\lambda_0}} \leq  A_2\sum_{n  = 1}^\infty n^{-2} < \infty
\]
if we choose $t_n =  c_1 \ln n$ with $c_1 \geq    2 /(\delta
\sqrt{\lambda_0})$. By the Borel-Cantelli lemma, there is a random
variable $N_1$ such that $(a')$ holds (with $N_1$ instead of $N$).

In order to establish $(b')$, note that for each $n$, the ball
$B((b-\delta) t_n)$ can be covered by the balls
$B^1_n,...,B^{m(n)}_n$ of radius $1/2$ centered at $t_n
v^1,...,t_n v^{m(n)}$ in such a way that the centers of the balls
are inside $B((b-\delta) t_n)$ and $m(n) = O(t_n^d)$. Let $g^k_n =
g(B(1/2), t_n, v^k)$. As in (\ref{ecooi}), there is $\gamma > 0$
such that
\[
\mathrm{E}(\frac{n^{x}_{t_n}(B^k_n) }{g^k_n} - \xi^x)^2 \leq  A
e^{- \gamma t_n}.
\]
By the Chebyshev inequality, for each $k = 1,...,m(n)$ and each $K
> 0$,
\[
\mathrm{P}(n^{x}_{t_n}(B^k_n)  = 0,~~\xi^x \geq K) \leq
\mathrm{P}(n^{x}_{t_n}(B^k_n)  \leq K g^k_n/2,~~\xi^x \geq K) \leq
\frac{ 4 A e^{- \gamma t_n}}{K^2}.
\]
Therefore,
\[
\sum_{n=1}^\infty \mathrm{P}(n^{x}_{t_n}(B^k_n)  = 0~~~{\rm
for}~{\rm some}~k,~~\xi^x \geq K) \leq \sum_{n=1}^\infty A_2(K)
e^{- \gamma t_n} t_n^d.
\]
The series on the right hand side of this inequality converges if
we choose $t_n = c_2 \ln n$ with $c_2 \geq 2/\gamma$. Since $K >
0$ was arbitrary, by the Borel-Cantelli lemma, $(b')$ holds (with
$N_2$ instead of $N$). It remains to take $c = \max(c_1, c_2)$ and
then $N = \max(N_1, N_2)$. \qed

\section{Limiting distribution in the case of finitely many
particles} \label{finitely}

In this section we make a couple of remarks concerning the
distribution of the total number of particles on the event $\xi^x
= 0$.
\\
\\
{\bf Remark.} Let $N^x_\infty = \lim_{t \rightarrow \infty}
N^x_t$. This random variable is finite almost surely on the event
$\xi^x = 0$ (this can be proved similarly to the way it was done
for the corresponding statement in the critical case in \cite{K}).
In other words, with probability one, the total number of
particles either grows exponentially or tends to a finite limit.
The latter event has nonzero probability if and only if $d \geq
3$.
\\
\\
{\bf Remark.} Let $d \geq 3$ and define $M^n(x) = \mathrm{P}(
N^x_\infty = n)$, $n \geq 1$. The quantities $M^n(x)$ satisfy a
recursive system of partial differential equations. Namely,
\begin{equation} \label{ft1}
\frac{1}{2} \Delta M^1(x) = v(x),
\end{equation}
with the condition at infinity
\[
\lim_{|x| \rightarrow \infty} M^1(x) = 1.
\]
For $n \geq 2$, we have
\begin{equation} \label{ft2}
\frac{1}{2} \Delta M^n(x) = v(x) \sum_{k = 1}^{n-1} M^k(x) M^{n-k}
(x),
\end{equation}
\[
\lim_{|x| \rightarrow \infty} M^n(x) = 0.
\]
Equations (\ref{ft1}) and (\ref{ft2}) can be obtained by
considering the behavior of the initial particle on the time
interval $[0,\delta]$ such that $\delta \downarrow 0$, with the
left hand side accounting for the diffusive motion and the right
hand side for the branching.


\begin{thebibliography}{9}


%\bibitem{Alb} Albeverio S., Gesztesy F.,  H{\o}egh-Krohn R.,
%              Holden H., {\it Solvable models in quantum mechanics},
% AMS Chelsea Publishing, Providence, RI, 2005.

\bibitem{ABY} Albeverio, S., Bogachev L., Yarovaya E., {\it Asymptotics
of branching symmetric random walk on the lattice with a single
source},  C. R. Acad. Sci. Paris Sér. I Math.  326  (1998),  no.
8, 975--980.

\bibitem{ABY2} Albeverio, S., Bogachev, L., Yarovaya, E., {\it Branching random walk with a single source},
Communications in difference equations (Poznan, 1998),  9--19,
Gordon and Breach, Amsterdam, 2000.

%\bibitem{Aronson1} Aronson D.G., Weinberger H.F., {\it Non-linear
%diffusion in population genetics, combustion and nerve
%propagation}. Proceedings of the Tulane program in partial
%differential equations (Lecture Notes in Mathematics),
%Springer-Verlag.

%\bibitem{Aronson2} Aronson D.G., Weinberger H.F., {\it
%Multi-dimensional non-linear diffusion in population genetics}.
%Advances in Math., 1978, 30, 33-76.

\bibitem{AH} Asmussen S., Hering H., {\it Strong Limit Theorems
for general supercritical branching processes with applications to
branching diffusions}, Z. Wahrscheinlichkeitstheorie und Verw.
Gebiete, 36 no 3, pp 195-212 (1976).

\bibitem{BY1} Bogachev L., Yarovaya E., {\it A limit theorem for a supercritical branching
 random walk on $Z^d$ with a single source},
translation in  Russian Math. Surveys  53  (1998),  no. 5,
1086--1088.

%\bibitem{CaMo} Carmona, R. A.; Molchanov, S. A., {\it Stationary parabolic Anderson
%model and intermittency.}  Probab. Theory Related Fields  102
%(1995),  no. 4, 433--453.

\bibitem{CKMV1} Cranston M., Koralov L., Molchanov S., Vainberg B., {\it
Continuous model for homopolymers},  Journal of Functional
Analysis 256 (2009), no. 8, 2656--2696.

%\bibitem{DaFl} Dawson D., Fleischmann K., {\it A super-Brownian motion with a single point catalyst},
%Stochastic Process. Appl.  49  (1994),  no. 1, 3--40.

%\bibitem{DFL} Dawson, D., Fleischmann K., Le Gall J., {\it Super-Brownian motions in catalytic media},
% Branching processes (Varna, 1993),  122--134, Lecture Notes in Statist., 99, Springer, New York, 1995.

%\bibitem{DFLM} Dawson, D., Fleischmann, K., Li Y., Mueller C., {\it Singularity of super-Brownian
%local time at a point catalyst},  Ann. Probab.  23  (1995),  no.
%1, 37--55.

\bibitem{EHK} Englander J., Harris S. C.,  Kyprianou A.E., {\it Strong Law of Large numbers for branching
diffusions}, Annales de l'Institut Henri Poincaré (B) Probability
and Statistics, (2010), Vol.46, No. 1, 279-298.

\bibitem{Fi1} Fischer R.A., {\it On the dominance ratio}, PRS
Edinburgh, 42 (1922), 321-341.
\bibitem{Fi2} Fischer R.A., {\it The genetical theory of natural
selection}, Oxford University Press, 1930.
\bibitem{Fi3} Fischer R.A., {\it The distribution of gene ratios
for rare mutations} PRS Edinburgh, 50 (1930) 204-219.

%\bibitem{FL} Fleischmann K., Le Gall J., {\it A new approach to
%the single point catalytic super-Brownian motion},  Probab. Theory
%Related Fields  102  (1995),  no. 1, 63--82

%\bibitem{Frank} Frank-Kamenezkii D.A., {\it Diffusion and heat
%transfer in chemical kinetics}, New York: Plenum Press.

%\bibitem{FSh} Frechet  M., Shohat J.,
%{\it A Proof of the Generalized Second-Limit Theorem in the Theory
%of Probability}, Transactions of the American Mathematical Society
%Vol. 33, No. 2 (Apr., 1931), pp. 533-543.

%\bibitem{GK} Gartner J., Konig W., {\it The parabolic Anderson
%model. Interacting stochastic systems}, 153--179, Springer,
%Berlin, 2005.

%\bibitem{GKM} Gartner, J., Konig, W., Molchanov, S. A., {\it Almost sure
%asymptotics for the continuous parabolic Anderson model},  Probab.
%Theory Related Fields  118  (2000),  no. 4, 547--573.

%\bibitem{GM1} Gartner, J., Molchanov, S. A., {\it Parabolic
%problems for the Anderson model. I. Intermittency and related
%topics},  Comm. Math. Phys.  132  (1990),  no. 3, 613--655.

%\bibitem{GM2} Gartner, J., Molchanov, S. A., {\it Parabolic problems for the Anderson
%model. II. Second-order asymptotics and structure of high peaks},
%Probab. Theory Related Fields  111  (1998),  no. 1, 17--55.

\bibitem{Hal} Haldane J.B.S., {\it A mathematical theory of natural and
artificial selection, part V: Selection and mutation}, PCPS, 23
(1927), 838-844.

\bibitem{HaU} Hawkins D., Ulam S., {\it Theory of Multiplicative
Processes I}, Los Alamos Scientific Laboratory, LADS-265 (1944).

%\bibitem{KLMS} Konig W., Lacoin H., Morters P., Sidorova. N.,
%{\it A two cities theorem for the parabolic Anderson model},  Ann.
%Probab. 37  (2009),  no. 1, 347--392.

%\bibitem{Kr}  Krasnoselski M., {\it Positive Solutions of Operator
%Equations}, Groningen, P. Noordhoff [1964].


\bibitem{K} Koralov L., {\it Branching
diffusion in inhomogeneous media}, submitted for publication.

\bibitem{Se} Semenov N.N., {\it Chain Reactions}, (in Russian),
Goshimizdat, 1934.


%\bibitem{Va} Vainberg B., {\it On the Short-Wave Asymptotic Behavior of
%Solutions of Stationary Problems and the Asymptotic Behavior as $t
%\rightarrow \infty$ of Solutions of Non-Stationary Problems},
%Russian Math Surveys, 30:2, 1975, pp 1-58.


\bibitem{W} Watanabe S., {\it A limit theorem of branching
processes and continuous state branching processes}, J. Math.
Kyoto Univ. 8, 141-167 (1968).

\bibitem{WG} Watson H.W., Galton F., {\it On the probability of
the extinction of families}, J. Anthropol. Inst. Great Britain and
Ireland, 4 (1874), 138-144.

\bibitem{Ya} Yarovaya E., {Branching random walks in inhomogeneous
media} (in Russian), Moscow University Meh-Mat publication.
\end{thebibliography}
\end{document}